# ON THE NEW BASES ATTACHED TO WEYL GROUPS OF EXCEPTIONAL TYPE

G. Lusztig

## Introduction

**0.1.** For a finite group $\Gamma$ let $\mathrm{Irr}(\Gamma)$ be the set of isomorphism classes of irreducible $\mathbf{C}[\Gamma]$-modules. For any finite group $\Gamma$, we denote by $M(\Gamma)$ the set of pairs $(x,\sigma)$ (up to $\Gamma$-conjugacy) where $x \in \Gamma$ and $\sigma \in \mathrm{Irr}(Z_\Gamma(x))$ ($Z_\Gamma(x)$ is the centralizer of $x$ in $\Gamma$). Let $\mathbf{C}[M(\Gamma)]$ be the $\mathbf{C}$-vector space with basis $M(\Gamma)$.

Following [L79] we define a pairing $\{:\} : M(\Gamma) \times M(\Gamma) \to \mathbf{C}$ by

$$\{(x,\sigma),(y,\tau)\} = \sum_{g\in\Gamma; xgyg^{-1}=gyg^{-1}x} |Z(x)|^{-1}|Z(y)|^{-1}\mathrm{tr}(g^{-1}x^{-1}g,\tau)\mathrm{tr}(gyg^{-1},\sigma).$$

We define a linear map (non-abelian Fourier transform) $A_\Gamma : \mathbf{C}[M(\Gamma)] \to \mathbf{C}[M(\Gamma)]$ by

$$A_\Gamma(y,\tau) = \sum_{(x,\sigma)\in M(\Gamma)} \{(x,\sigma),(y,\tau)\}(x,\sigma).$$

Let

$$\mathbf{R}_{\geq 0}[M(\Gamma)] = \sum_{(x,\sigma)\in M(\Gamma)} \mathbf{R}_{\geq 0}(x,\sigma) \subset \mathbf{C}[M(\Gamma)].$$

We say that $f \in \mathbf{C}[M(\Gamma)]$ is *bipositive* if $f \in \mathbf{R}_{\geq 0}[M(\Gamma)]$ and $A_\Gamma(f) \in \mathbf{R}_{\geq 0}[M(\Gamma)]$.

Let $\tilde{Z}(\Gamma)$ be the set of all pairs $(H,H')$ of subgroups of $\Gamma$ (up to $\Gamma$-conjugacy) where $H$ is a normal subgroup of $H'$.

For $(H,H') \in \tilde{Z}(\Gamma)$ let $\mathbf{s}_{H,H';\Gamma} : \mathbf{C}[M(H'/H)] \to \mathbf{C}[M(\Gamma)]$ be the linear map introduced in [L20, 3.1]. It is known that

(a) $$A_\Gamma \mathbf{s}_{H,H';\Gamma} = \mathbf{s}_{H,H';\Gamma} A_{H'/H}.$$

It follows that

(b) $\mathbf{s}_{H,H';\Gamma}$ applied to any bipositive element of $\mathbf{C}[M(H'/H)]$ is a bipositive element of $\mathbf{C}[M(\Gamma)]$.

**0.2.** The main result of this paper is a construction of a $\mathbf{C}$-basis $B(\Gamma)$ of $\mathbf{C}[M(\Gamma)]$ in the case where $\Gamma$ is one of the symmetric groups $S_n, n \in \{1,2,3,4,5\}$ with the following properties:

(I) Any $b \in B(\Gamma)$ is bipositive.

(II) There is a unique bijection $\epsilon : M(\Gamma) \xrightarrow{\sim} B(\Gamma)$ such that for any $(x,\sigma) \in M(\Gamma)$, $(x,\sigma)$ appears with nonzero coefficient in $\epsilon(x,\sigma)$. (That coefficient is actually 1.)

(III) The matrix expressing the basis $B(\Gamma)$ in terms of the obvious basis of $\mathbf{C}[M(\Gamma)]$ is upper triangular (for some order on $M(\Gamma)$) with diagonal entries 1.

(IV) The matrix of $A_\Gamma : \mathbf{C}[M(\Gamma)] \to \mathbf{C}[M(\Gamma)]$ with respect to the basis $B(\Gamma)$ is upper triangular (for another order on $M(\Gamma)$) with diagonal entries $\pm 1$.

Our definition of $B(\Gamma)$ has two main ingredients. The first one is a certain subset $Z(\Gamma)$ of $\tilde{Z}(\Gamma)$. The set $Z(\Gamma)$ (without the pair $(1,1)$ when $n \in \{4,5\}$) has already been introduced in [L19] where it was used to define a part of $B(\Gamma)$; a uniform description of $Z(\Gamma)$ was given in [L24].

The second ingredient is a definition of a small subset $Prim(H'/H)$ (said to be the primitive elements) of $\mathbf{C}[M(H'/H)]$ for $(H,H') \in Z_\Gamma$. Then we define $B(\Gamma)$ to be the subset of $\mathbf{C}[M(\Gamma)]$ consisting of all $\mathbf{s}_{H,H';\Gamma}P$ where $(H,H')$ runs through $Z(\Gamma)$ and $P$ runs through $Prim(H'/H)$. The proof that $B(\Gamma)$ is a basis with properties (I)-(IV) is given in 3.8, 3.10.

**0.3.** When $n$ is $3,4,5$ a basis of $\mathbf{C}[M(\Gamma)]$ was given in [L20] (but for $n=5$, [L20] contained an error). The basis defined in this paper differs from that in [L20], although several of the ideas in [L20] are still used. This paper can be regarded as a substitute for [L20, 3.2-3.9] and for [L20a,§5].

## 1. Cyclic groups and Fourier transform

**1.1.** Let $n \in \{2,3,4,\dots\}$. Let $C = \mathbf{Z}/n\mathbf{Z}$ with the group structure being the addition. Let $\tilde{C}$ be the set of generators of this group. It is a group under multiplication. We have an isomorphism $C \to \mathrm{Hom}(C, \mathbf{C}^*)$ given by $a \mapsto [z \mapsto \xi_n^{az}]$ where $\xi_n = e^{2\pi i/n}$. Hence we can identify $M(C) = C \times C$.

Let $C' = C - \{0\}$ and let

$$\mathcal{S} = \sum_{(k,l)\in C'\times C'} (k,l) + (n-1)(0,0) \in \mathbf{C}[M(C)].$$

We show:

(a) $A_C(\mathcal{S}) = \mathcal{S}$.

We have

$$
\begin{aligned}
nA_C(\mathcal{S}) &= \sum_{(k,l)\in C'\times C',(k',l')\in C\times C} \xi_n^{-kl'+k'l}(k',l') + (n-1)\sum_{(k',l')\in C\times C}(k',l') \\
&= \sum_{(k',l')\in C'\times C'}\sum_{(k,l)\in C'\times C'} \xi_n^{-kl'+k'l}(k',l') + \sum_{l'\in C'}\sum_{(k,l)\in C'\times C'}\xi_n^{-kl'}(0,l') \\
&+ \sum_{k'\in C'}\sum_{(k,l)\in C'\times C'}\xi_n^{k'l}(k',0) + \sum_{(k,l)\in C'\times C'}(0,0) + (n-1)\sum_{(k',l')\in C\times C}(k',l') \\
&= \sum_{(k',l')\in C'\times C'}(-1)(-1)(k',l') + \sum_{l'\in C'}(n-1)(-1)(0,l') \\
&+ \sum_{k'\in C'}(n-1)(-1)(k',0) + (n-1)^2(0,0) + (n-1)\sum_{(k',l')\in C\times C}(k',l') \\
&= \sum_{(k',l')\in C'\times C'}(1+n-1)(k',l') + ((n-1)^2+n-1)(1,1) = n\mathcal{S}.
\end{aligned}
$$

This proves (a).

For any $a\in\tilde{C}$ let

$$\mathcal{S}_a = \sum_{(k,l)\in\tilde{C}\times\tilde{C});kl=a}(k,l)\in\mathbf{C}[M(C)].$$

We have

$$nA_C(\mathcal{S}_a) = \sum_{(k',l')\in C\times C} c_{k',l'}(k',l')$$

where

$$c_{k',l'} = \sum_{(k,l)\in\tilde{C}\times\tilde{C};kl=a}\xi_n^{-kl'+k'l}.$$

If $k'l'=b\in\tilde{C}$ then (setting $u=k/k'$) we have

(b) $$c_{k',l'} = \sum_{u\in\tilde{C}}\xi_n^{-ub+au^{-1}}.$$

We set $\mathcal{T}=\sum_{k\in\tilde{C}}(k,0)$, $\mathcal{T}'=\sum_{l\in\tilde{C}}(0,l)$.

For example if $n=2$, we have $\xi_2=-1$, $C'=\tilde{C}=\{1\}$ and

(c) $$A_C(\mathcal{S}_1+(0,0)) = \mathcal{S}_1+(0,0).$$

**1.2.** We now assume that $n=3$. Let $\theta=\xi_3$. We have $C'=\tilde{C}=\{1,2\}$. We show

(a) $$A_C(\mathcal{S}_1-\mathcal{S}_2)=\mathcal{S}_1-\mathcal{S}_2.$$

Indeed, using 1.1(b) we have

$$3A_C(\mathcal{S}_1-\mathcal{S}_2)=\sum_{b\in\tilde{C}}(\sum_{u\in\tilde{C}}(\theta^{-ub+u^{-1}}-\theta^{-ub+2u^{-1}})\mathcal{S}_b$$
$$=\sum_{b\in\tilde{C}}(\theta^{-b+1}-\theta^{-b+2}+\theta^{b-1}-\theta^{b-2})\mathcal{S}_b=(\theta^{1-1}-\theta^{-1+2}+\theta^{-1+1}-\theta^{1-2})\mathcal{S}_1+$$
$$(\theta^{1+1}-\theta^{1+2}+\theta^{-1-1}-\theta^{-1-2})\mathcal{S}_{-1}$$
$$=(1-\theta^{-1}+1-\theta)\mathcal{S}_1+(\theta-1+\theta^{-1}-1)\mathcal{S}_{-1}=3(\mathcal{S}_1-\mathcal{S}_2),$$

as required. From 1.1(a) we have

(b) $$A_C(\mathcal{S}_1+\mathcal{S}_2+2(0,0))=\mathcal{S}_1+\mathcal{S}_2+2(0,0).$$

This, together with (a) implies

(c) $$A_C(\mathcal{S}_1+(0,0))=\mathcal{S}_1+(0,0),$$

(d) $$A_C(\mathcal{S}_2+(0,0))=\mathcal{S}_2+(0,0).$$

Since

$$3A_C(0,0)=\sum_{(k,l)\in C\times C}(k,l),$$

we see that

(e) $A_C(\mathcal{S}_1+2(0,0))=\mathcal{S}_1+2(0,0) \mod \mathbf{R}(\mathcal{S}_1+\mathcal{S}_2+2(0,0))+\mathbf{R}(\mathcal{T}+\mathcal{T}')+\mathbf{R}(0,0).$

**1.3.** We now assume that $n=4$. We have $\xi_4=i$. We have $\tilde{C}=\{1,-1\}$. We show:

(a) $$A_C(\mathcal{S}_1-\mathcal{S}_{-1})=\mathcal{S}_1-\mathcal{S}_{-1}.$$

Indeed,

$$4A_C(\mathcal{S}_1-\mathcal{S}_{-1})=\sum_{b\in\tilde{C}}(\sum_{u\in\tilde{C}}(i^{-ub+u^{-1}}-i^{-ub+3u^{-1}})\mathcal{S}_b$$
$$=\sum_{b\in\tilde{C}}(i^{-b+1}-i^{-b+3}+i^{b-1}-i^{b-3})\mathcal{S}_b$$
$$=(i^{1-1}-i^{-1+3}+i^{-1+1}-i^{1-3})\mathcal{S}_1$$
$$+(i^{1+1}-i^{1+3}+i^{-1-1}-i^{-1-3})\mathcal{S}_{-1}=4(\mathcal{S}_1-cs_{-1}),$$

as required.

We set $\mathcal{S}'=\mathcal{S}_1+\mathcal{S}_{-1}+2(2,2)+2(0,0)$. We show

(b) $$A_C(\mathcal{S}')=\mathcal{S}'.$$

We have

$$4A_C(\mathcal{S}')=\sum_{(k,l)\in\tilde C\times\tilde C,(k',l')\in C\times C} i^{-kl'+k'l}(k',l')+\sum_{(k',l')\in C\times C}(i^{2l'+2k'}+1)(k',l')$$
$$=\sum_{(k',l')\in\tilde C\times\tilde C}\sum_{(k,l)\in\tilde C\times\tilde C} i^{-kl'+k'l}(k',l')+\sum_{l'\in\tilde C}\sum_{(k,l)\in\tilde C\times\tilde C} i^{-kl'+2l}(2,l')$$
$$+\sum_{l'\in\tilde C}\sum_{(k,l)\in\tilde C\times\tilde C} i^{-kl'}(0,l')+\sum_{k'\in\tilde C}\sum_{(k,l)\in\tilde C\times\tilde C} i^{-2k+k'l}(k',2)$$
$$+\sum_{k'\in\tilde C}\sum_{(k,l)\in\tilde C\times\tilde C} i^{k'l}(k',0)+\sum_{(k,l)\in\tilde C\times\tilde C} i^{-2k+2l}(2,2)$$
$$+\sum_{(k,l)\in\tilde C\times\tilde C} i^{2l}(2,0)+\sum_{(k,l)\in\tilde C\times\tilde C} i^{-2k}(0,2)$$
$$+\sum_{(k,l)\in\tilde C\times\tilde C}(0,0)+\sum_{(k',l')\in\tilde C\times\tilde C} i^{-2l'+2k})(k',l')$$
$$+2\sum_{(k',l')\in\tilde C\times\tilde C}(k',l')+2(i^{2.2+2.2}+1)(2,2)$$
$$+2(i^{2.2}+1)(2,0)+(i^{2.2}+1)(0,2)+2(1+1)(0,0)$$
$$=0+0+0+0+0+4(2,2)-4(2,0)-4(0,2)+4(0,0)$$
$$+2\sum_{(k',l')\in\tilde C\times\tilde C}(k',l')++2\sum_{(k',l')\in\tilde C\times\tilde C}(k',l')$$

(b)
$$+4(2,2)+4(2,0)+4(0,2)+4(0,0)=4\mathcal{S}',$$

as required.

From (a),(b) we deduce

(c) $$A_C(\mathcal{S}_1+(2,2)+(0,0))=\mathcal{S}_1+(2,2)+(0,0),$$

(d) $$A_C(\mathcal{S}_{-1}+(2,2)+(0,0))=\mathcal{S}_{-1}+(2,2)+(0,0).$$

We have

$$4A_C(0,0)=\sum_{(k,l)\in C\times C}(k,l)$$
$$=\mathcal{S}_1+\mathcal{S}_{-1}+(0,0)+(0,2)+(2,0)+(2,2)+\sum_{(k,l)\in C\times C;k+l=1 \mod 2}(k,l),$$

$$\begin{aligned}4A_C(2,2) &= \sum_{(k,l)\in C\times C} i^{2k+2l}(k,l)\\ &= \sum_{(k,l)\in C\times C; k+l=0 \mod 2}(k,l) - \sum_{(k,l)\in C\times C; k+l=1 \mod 2}(k,l)\\ &= \mathcal{S}_1+\mathcal{S}_{-1}+(0,0)+(0,2)+(2,0)+(2,2)-\sum_{(k,l)\in C\times C; k+l=1 \mod 2}(k,l),\end{aligned}$$

hence

$$4A_C((2,2)+(0,0)) = 2(\mathcal{S}_1+\mathcal{S}_{-1}+(0,0)+(0,2)+(2,0)+(2,2)).$$

We see that

$$\begin{aligned}&A_C(\mathcal{S}_1+2(2,2)+2(0,0)) = \mathcal{S}_1+2(2,2)+2(0,0)\\ &\mod \mathbf{R}(\mathcal{S}_1+\mathcal{S}_{-1}+2(2,2)+2(0,0))+\mathbf{R}(0,0)+\mathbf{R}(0,2)+\mathbf{R}(2,0)+\mathbf{R}(2,2).\end{aligned} \tag{e}$$

**1.4.** We now assume that $n=5$. Let $\zeta=\xi_5$. We have $C'=\tilde{C}=\{1,2,3,4\}$. For $k\in\mathbf{Z}$ we set $[k]=\zeta^k+\zeta^{-k}$. We have $[1]+[2]=-1, [1]-[2]=\sqrt{5}>0$. In particular, $[1]=(\sqrt{5}-1)/2>0, -[2]=(1+\sqrt{5})/2>0$. We have

$$5A_C(0,0)=\sum_{a\in\tilde{C}}\mathcal{S}_a+\mathcal{T}+\mathcal{T}'+(0,0). \tag{a}$$

For $a\in\tilde{C}$ we have (using 1.1(b)):

$$\begin{aligned}5A_C(\mathcal{S}_a) &= \sum_{b\in\tilde{C}}(\sum_{u\in\tilde{C}}\zeta^{-ub+au^{-1}})\mathcal{S}_b-\mathcal{T}-\mathcal{T}'+4(0,0)\\ &= \sum_{b\in\tilde{C}}([b-a]+[2b-3a])\mathcal{S}_b-\mathcal{T}-\mathcal{T}'+4(0,0)\\ &= ([1-a]+[2-3a])\mathcal{S}_1+([2-a]+[4-3a])\mathcal{S}_2+([3-a]+[1-3a])\mathcal{S}_3\\ &\quad+([4-a]+[3-3a])\mathcal{S}_4-\mathcal{T}-\mathcal{T}'+4(0,0).\end{aligned} \tag{b}$$

From this and (a) we deduce

$$A_C(\mathcal{S}_2+\mathcal{S}_3+2(0,0))=\mathcal{S}_2+\mathcal{S}_3+2(0,0), \tag{c}$$

$$A_C(\mathcal{S}_1+\mathcal{S}_4+2(0,0))=\mathcal{S}_1+\mathcal{S}_4+2(0,0). \tag{d}$$

It follows that

(e) $$A_C(\mathcal{S}_1+\mathcal{S}_2+\mathcal{S}_3+\mathcal{S}_4+4(0,0))=\mathcal{S}_1+\mathcal{S}_2+\mathcal{S}_3+\mathcal{S}_4+4(0,0),$$

which also follows from 1.1(a).

We also see that

(f) $$A_C(\mathcal{S}_2+\mathcal{S}_3+4(0,0))= \mathcal{S}_2+\mathcal{S}_3+4(0,0) \mod \mathbf{R}(\mathcal{S}_1+\mathcal{S}_2+\mathcal{S}_3+\mathcal{S}_4+4(0,0))+\mathbf{R}(\mathcal{T}+\mathcal{T}')+\mathbf{R}(0,0).$$

From (a),(b) we have

$$\begin{aligned}&5A_C(-[2]\mathcal{S}_1+\mathcal{S}_2+[2]^2(0,0))\\ &=(-[2][0]-[2][1]+[1]+[1]+[2]^2)\mathcal{S}_1+(-[2][1]-[2][1]+[0]+[2]+[2]^2)\mathcal{S}_2\\ &+(-[2]^2-[2]^2+[1]+[0]+[2]^2)\mathcal{S}_3+(-[2]^2-[2][0]+[2]+[2]+[2]^2)\mathcal{S}_4\\ &+([2]^2-(-[2]+1))\mathcal{T}+([2]^2-(-[2]+1))\mathcal{T}'+([2]^2+(-[2]+1)4)(0,0).\end{aligned}$$

Hence

(g) $$A_C(-[2]\mathcal{S}_1+\mathcal{S}_2+[2]^2(0,0)=-[2]\mathcal{S}_1+\mathcal{S}_2+[2]^2(0,0).$$

We see that

(h) $$\begin{aligned}&A_C(-[2]\mathcal{S}_1+\mathcal{S}_2+4(0,0))\\ &=(-[2]\mathcal{S}_1+\mathcal{S}_2+4(0,0)) \mod \mathbf{R}(\mathcal{S}_1+\mathcal{S}_2+\mathcal{S}_3+\mathcal{S}_4+4(0,0))\\ &+\mathbf{R}(\mathcal{T}+\mathcal{T}')+\mathbf{R}(0,0).\end{aligned}$$

We show that

$$-[2]\mathcal{S}_1+\mathcal{S}_2+4(0,0)\in\mathbf{C}[M(C)] \text{ is bipositive.}$$

Since $-[2]>0$ it is enough to show that

$$A_C(-[2]\mathcal{S}_1+\mathcal{S}_2+4(0,0))\in\mathbf{R}_{\ge 0}[M(C)]$$

or that

$$-[2]\mathcal{S}_1+\mathcal{S}_2+[2]^2(0,0)+A_C((4-[2]^2)(0,0))\in\mathbf{R}_{\ge 0}[M(C)].$$

(To see this we use that $4-[2]^2>0$.)

From (a) we have

$$5A_C(\mathcal{S}_1)=R\mathcal{S}_1+M(-[2]\mathcal{S}_1+\mathcal{S}_2)+T(\mathcal{S}_1+\mathcal{S}_2+\mathcal{S}_3+\mathcal{S}_4)+U(\mathcal{S}_2+\mathcal{S}_3)-\mathcal{T}-\mathcal{T}'+4(0,0)$$

where

$$R-[2]M+T=[1-1]+[2-3]=[0]+[1],$$
$$M+U+T=[1]+[1],T=[2]+[0],U+T=[2]+[2]$$

so that $M=2[1]-2[2]$,

$$R=[0]+[1]-[2]-[0]+2[1][2]-2[2][2]$$
$$=[1]-[2]-2-2([1]+[0])=-[1]-[2]-6=1-6=-5.$$

We see that

$$A_C(\mathcal{S}_1+4(0,0))=-(\mathcal{S}_1+4(0,0))$$
$$\text{mod } \mathbf{R}(-[2]\mathcal{S}_1+\mathcal{S}_2+4(0,0))+\mathbf{R}(\mathcal{S}_2+\mathcal{S}_3+4(0,0))$$
$$+\mathbf{R}(\mathcal{S}_1+\mathcal{S}_2+\mathcal{S}_3+\mathcal{S}_4+4(0,0))+\mathbf{R}(\mathcal{T}+\mathcal{T}')+\mathbf{R}(0,0). \tag{j}$$

We show:

$$\mathcal{S}_1+4(0,0)\in\mathbf{C}[M(C)] \text{ is bipositive.}$$

It is enough to show that $A_C(\mathcal{S}_1+4(0,0))\in\mathbf{R}_{\ge0}[M(C)]$ or that

$$([1-1]+[2-3])\mathcal{S}_1+([2-1]+[4-3])\mathcal{S}_2+([3-1]+[1-3])\mathcal{S}_3$$
$$+([4-1]+[3-3])\mathcal{S}_4-\mathcal{T}-\mathcal{T}'+4(0,0)+4(\mathcal{S}_1+\mathcal{S}_2+\mathcal{S}_3+\mathcal{S}_4+\mathcal{T}$$
$$+\mathcal{T}'+(0,0))\in\mathbf{R}_{\ge0}[M(C)].$$

This follows from the fact that $[k]+[l]+4\ge0$ for any $k,l$.

## 2. Primitive elements of $M(\Gamma)$

**2.1.** In this section $\Gamma$ denotes one of the groups $S_1,S_2,S_3,S_4,S_5,S_3\times S_2,S_2\times S_2$. We shall define a subset $Prim(\Gamma)$ of $\mathbf{C}[M(\Gamma)]$ and a function $P\mapsto \text{sg}'(P)$ with values in $\{\pm1\}$ on it. (We shall use the notation of [L84] for the elements of $M(S_n)$.)

We set $Prim(S_1)=\{P_1\}$ where $P_1=(1,1)$, $\text{sg}'=1$.

We set $Prim(S_2)=\{P_{-1},P_{(-1)^0}\}$, where
$P_{-1}=(g_2,\epsilon)+(1,1)$, $\text{sg}'=1$,
$P_{(-1)^0}=(1,1)$, $\text{sg}'=-1$.

We set $Prim(S_3) = \{P_{\theta^2}, P_{\theta^0}, P_\theta\}$ where
$P_{\theta^2} = (g_3, \theta) + (g_3, \theta^2) + \mathcal{H}_3$, $\mathrm{sg}' = 1$,
$P_{\theta^0} = (1, 1)$, $\mathrm{sg}' = -1$,
$P_\theta = (g_3, \theta) + \mathcal{H}'_3$, $\mathrm{sg}' = 1$,
and $\mathcal{H}_3 = (1, 1) + (1, \epsilon)$, $\mathcal{H}'_3 = (1, 1)$.

We set $Prim(S_4) = \{P_{i^3}, P_{i^0}, P_i\}$ where
$P_{i^3} = (g_4, i) + (g_4, i^3) + \mathcal{H}_4$, $\mathrm{sg}' = 1$,
$P_{i^0} = (1, 1)$, $\mathrm{sg}' = 1$,
$P_i = (g_4, i) + \mathcal{H}'_4$, $\mathrm{sg}' = 1$,
and $\mathcal{H}_4 = (g'_2, \epsilon')+(g'_2, \epsilon'')+(1, \lambda^3)+(1, \sigma)+(1, 1)$, $\mathcal{H}'_4 = (g'_2, \epsilon')+(1, \sigma)+(1, 1)$.

We set $Prim(S_5) = \{P_{\zeta^4}, P_{\zeta^0}, P_{\zeta^3}, P_{\zeta^2}, P_\zeta\}$ where
$P_{\zeta^4} = (g_5, \zeta) + (g_5, \zeta^2) + (g_5, \zeta^3) + (g_5, \zeta^4) + \mathcal{H}_5$, $\mathrm{sg}' = 1$,
$P_{\zeta^0} = (1, 1)$, $\mathrm{sg}' = -1$,
$P_{\zeta^3} = (g_5, \zeta^2) + (g_5, \zeta^3) + \mathcal{H}'_5$, $\mathrm{sg}' = 1$,
$P_{\zeta^2} = -(\zeta^2 + \zeta^3)(g_5, \zeta) + (g_5, \zeta^2) + \mathcal{H}_5$, $\mathrm{sg}' = 1$,
$P_\zeta = (g_5, \zeta) + \mathcal{H}_5$, $\mathrm{sg}' = -1$,
and $\mathcal{H}_5 = (1, 1)+(1, \nu)+2(1, \lambda^2)+(1, \nu')+(1, \lambda^4)$, $\mathcal{H}'_5 = (1, 1)+(1, \nu)+(1, \lambda^2)$.

We set

$$Prim_{S_3\times S_2} = \{P_{(-\theta)^0}, P_{(-\theta)^1}, P_{(-\theta)^2}, P_{(-\theta)^3}, P_{(-\theta)^4}, P_{(-\theta)^5}\}$$

where

$$P_{(-\theta)^k} = P_{\theta^k} \otimes P_{(-1)^k} \in \mathbf{C}[M(S_3)] \otimes \mathbf{C}[M(S_2)] = \mathbf{C}[M(S_3 \times S_2)]$$

for $k = 0, 1, 2, 3, 4, 5$ and the values of $\mathrm{sg}'$ are $1, 1, -1, -1, -1, 1$ respectively.

We set

$$Prim_{S_2\times S_2} = \{P_{-1,-1}, P_*, P_{1,1}\}$$

where

$$P_{-1,-1} = (g_2, \epsilon) \otimes (g_2, \epsilon) + (g_2, \epsilon) \otimes (1, 1) + (1, 1) \otimes (g_2, \epsilon) + (1, 1) \otimes (1, 1),$$

$$P_* = (g_2, \epsilon) \otimes (g_2, \epsilon) + (1, 1) \otimes (1, 1),$$

$$P_{1,1} = (1, 1) \otimes (1, 1)$$

are in $\mathbf{C}[M(S_2)]\otimes\mathbf{C}[M(S_2)] = \mathbf{C}[M(S_2 \times S_2)]$ and the values of $\mathrm{sg}'$ are $1, -1, -1$ respectively.

**2.2.** For $n\in\{2,3,4,5\}$ we write $C_n$ instead of $C$ in §1.

If $n=2$ we have $C_2=S_2$. From the definitions we have

(a) $$\mathbf{s}_{(1,C);S_2}(\mathcal{S}_1+(0,0))=(g_2,\epsilon)+(1,1).$$

If $n=3$ we can regard $C_3$ as a subgroup of $S_3$. From the definitions we have

$$\mathbf{s}_{(1,C_2);S_3}(\mathcal{S}/2)=(g_3,\theta)+(g_3,\theta^2)+\mathcal{H}_3=P_{\theta^2},$$

$$\mathbf{s}_{(1,C_2);S_3}((\mathcal{S}_1+2(0,0))/2)=(g_3,\theta)+\mathcal{H}_3=P_\theta+\mathcal{H}_3-\mathcal{H}'_3.$$

We can regard $C_4$ as a subgroup of $S_4$. From the definitions we have

$$\mathbf{s}_{(1,C_4);S_4}(\mathcal{S}'/2)=(g_4,i)+(g_4,i^3)+\mathcal{H}_4=P_{i^3},$$

$$\mathbf{s}_{(1,C_4);S_4}((\mathcal{S}_1+2(2,2)+2(0,0))/2)=(g_4,i)+\mathcal{H}_4=P_i+\mathcal{H}_4-\mathcal{H}'_4.$$

We can regard $C_5$ as a subgroup of $S_5$. From the definitions we have

$$\mathbf{s}_{(1,C_5);S_5}(\mathcal{S}/4)=(g_5,\zeta)+(g_5,\zeta^2)+(g_5,\zeta^3)+(g_5,\zeta^4)+\mathcal{H}_5=P_{\zeta^4},$$

$$\mathbf{s}_{(1,C_5);S_5}((\mathcal{S}_2+\mathcal{S}_3+4(0,0))/4)=(g_5,\zeta^2)+(g_5,\zeta^3)+\mathcal{H}_5=P_{\zeta^3}+\mathcal{H}_5-\mathcal{H}'_5,$$

$$\mathbf{s}_{(1,C_5);S_5}((-[2]\mathcal{S}_1+\mathcal{S}_2+4(0,0))/4)=[-2](g_5,\zeta)+(g_5,\zeta^2)+\mathcal{H}_5=P_{\zeta^2},$$

$$\mathbf{s}_{(1,C_5);S_5}((\mathcal{S}_1+4(0,0))/4)=(g_5,\zeta)+\mathcal{H}_5=P_\zeta.$$

**Proposition 2.3.** *If $P\in PrimM(\Gamma)$ then $P$ is bipositive.*

We have

(a) $$A_\Gamma(1,1)=\sum_{(x,\sigma)\in M(\Gamma}|Z_\Gamma(x)|^{-1}\dim(\sigma)(x,\sigma).$$

We see that the proposition holds when $P=(1,1)$. We now assume that $P\neq(1,1)$. Assume first that $\Gamma=S_n$ with $n\in\{2,3,4,5\}$ and $P$ is one of

$P_{-1},P_{\theta^2},P_{i^3},P_{\zeta^4},P_\zeta,P_{\zeta^2}$.

Using the description of $P$ in terms of $C_n$ (in 2.2) and the fact that $\mathbf{s}_{1,C_n;S_n}$ takes bipositive elements of $M(C_n)$ to bipositive elements in $M(S_n)$ (see 0.1(b)) we are reduced to the bipositivity statements in §1. If $P$ is one of $P_\theta,P_i,P_{\zeta^3}$, the desired result is proved by direct computation. The case where $\Gamma=S_3\times S_2$ can be reduced to the cases of $S_3$ and $S_2$. If $\Gamma=S_2\times S_2$ and $P=P_{-1,-1}$ then $P$ is fixed by $A_\Gamma(P)$ hence is bipositive; if $P=P_*$ then all coefficients of $A_\Gamma(P)$ are in $\{0,1/2\}$ hence $P$ is bipositive.

**2.4.** For $\Gamma$ as in 2.1 we define a permutation $P \mapsto P^!$ of $Prim(\Gamma)$ as follows:
$P_1 \mapsto P_1$ if $\Gamma = S_1$
$P_{-1} \mapsto P_{(-1)^0} \mapsto P_{-1}$ if $\Gamma = S_2$
$P_{\theta^2} \mapsto P_\theta \mapsto P_{\theta^0} \mapsto P_{\theta^2}$ if $\Gamma = S_3$
$P_i \mapsto P_{i^0} \mapsto P_{i^3} \mapsto P_i$ if $\Gamma = S_4$
$P_{\zeta^4} \mapsto P_\zeta \mapsto P_{\zeta^0} \mapsto P_{\zeta^4}$,
$P_{\zeta^3} \mapsto P_{\zeta^2} \mapsto P_{\zeta^3}$ if $\Gamma = S_5$
$P_{-1,-1} \mapsto P_{1,1} \mapsto P_{-1,-1}$, $P_* \mapsto P_*$ if $\Gamma = S_2 \times S_2$

$$P_{(-\theta)^5} \mapsto P_{(-\theta)^4} \mapsto P_{(-\theta)^3} \mapsto P_{(-\theta)^2} \mapsto P_{(-\theta)^1} \mapsto P_{(-\theta)^0} \mapsto P_{(-\theta)^5}$$

if $\Gamma = S_3 \times S_2$.

## 3. The basis $B(\Gamma)$

**3.1.** If $\Gamma = S_n, n \in \{1, 2, 3, 4, 5\}$, we associate to $\Gamma$ the following set $Z(S_n)$ of pairs $(H, H')$ of subgroups of $\Gamma$.
$Z(S_1)$ is: $(1, 1)$.
$Z(S_2)$ is:

$$(1, S_2), (S_2, S_2), (1, 1).$$

$Z(S_3)$ is:

$$(1, S_3), (1, S_2), (S_2, S_2), (S_3, S_3), (1, 1).$$

Here $S_2$ is viewed as the subgroup of $S_3$ generated by a transposition.
$Z(S_4)$ is:

$$(1, S_4), (1, S_3), (1, S_2 \times S_2), (1, S_2), (S_2, S_2 \times S_2)$$
$$, (S_2 \times S_2, \Delta), (S_2, S_2), (S_3, S_3), (S_4, S_4), (S_2 \times S_2, S_2 \times S_2),$$
$$(\Delta, \Delta), (1, 1).$$

Here $S_3$ is identified with a subgroup of $S_4$, $S_2 \times S_2$ is viewed as the subgroup of $S_4$ generated by two distinct commuting transpositions $t, t'$, $\Delta$ is the centralizer of $tt'$ and $S_2$ is the subgroup generated by $t$.
$Z(S_5)$ is:

$$(1, S_5), (1, S_4), (1, S_3 \times S_2), (1, S_3), (1, S_2 \times S_2), (1, S_2),$$
$$(S_2, S_3 \times S_2), (S_2, S_2 \times S_2), (S_3, S_3 \times S_2), (S_2 \times S_2, \Delta),$$
$$(S_2 \times S_2, S_2 \times S_2), (S_3 \times S_2, S_3 \times S_2), (S_2, S_2), (S_3, S_3),$$
$$(S_4, S_4), (S_5, S_5), (\Delta, \Delta), (1, 1).$$

Here $S_4$ is identified with a subgroup of $S_5$, $S_3$ is identified with a subgroup of $S_5$, $S_2$ is the subgroup of $S_5$ generated by a transposition $t$ not in $S_3$ but commuting with $S_3$ so that $S_3 \times S_2$ is a subgroup of $S_5$; $S_2 \times S_2$ is identified with the subgroup generated by $t$ and by a transposition $t' \neq t$ commuting with $t$; $\Delta$ is the centralizer of $tt'$ in $S_5$.

If $\Gamma = S_2 \times S_2$ we associate to $\Gamma$ the set $Z(S_2 \times S_2)$ of pairs $(H, H')$ of subgroups of $\Gamma$ consisting of

$$(1, S_2 \times S_2), (S_2 \times 1, S_2 \times S_2), (1 \times S_2, S_2 \times S_2),$$
$$(1, S_2 \times 1), (1, 1 \times S_2), (S_2 \times S_2, S_2 \times S_2)$$
$$(S_2 \times 1, S_2 \times 1), (1 \times S_2, 1 \times S_2), (D, D), (1, 1)$$

where $D$ is the diagonal in $S_2 \times S_2$.

If $\Gamma = S_3 \times S_2$ we associate to $\Gamma$ the set $Z(S_3 \times S_2)$ of pairs $(H, H')$ of subgroups of $\Gamma$ consisting of $(H_1 \times H_2, H_1' \times H_2')$ where $(H_1, H_1') \in Z(S_3)$, $(H_2, H_2') \in Z(S_2)$.

If $\Gamma$ is one of

(a) $S_1, S_2, S_3, S_4, S_5, S_3 \times S_2, S_2 \times S_2$,

then for $(H, H') \in Z(\Gamma)$, $H'/H$ is one of the groups in (a), hence $Prim(H'/H) \subset \mathbf{C}[M(H'/H)]$ is defined.

Let $B'(\Gamma)$ be the set of all triples $(H, H', P)$ where $(H, H') \in Z(\Gamma)$ and $P \in Prim(H'/H)$.

We now assume that $\Gamma = S_n$ where $n \in \{1, 2, 3, 4, 5\}$. Let $Y_n$ be the totally ordered set consisting of

(1) if $n = 1$

$(1) < (-1)$ if $n = 2$

$(1) < (-1) < (\theta) < (\theta^2)$ if $n = 3$

$(1) < (-1) < (-1') < (\theta) < (\theta^2) < (i) < (i^3)$ if $n = 4$

$(1) < (-1) < (-1') < (\theta) < (\theta^2) < (i) < (i^3) < (\zeta) < (\zeta^2) < (\zeta^3) < (\zeta^4) < (-\theta) < (-\theta^2)$ if $n = 5$.

We define $\mathbf{c} : B'(S_n) \to Y_n$ by $\mathbf{c}(H, H', P) = (\xi)$ if $H'/H \neq S_2 \times S_2$, $P = P_\xi$, while if $H'/H = S_2 \times S_2$, we have

$\mathbf{c}(H, H', P_{1,1}) = (1), \mathbf{c}(H, H', P_{-1,-1}) = (-1'), \mathbf{c}(H, H', P_*) = (-1)$.

For $(H, H', P), (\tilde{H}, \tilde{H}', \tilde{P})$ in $B'(\Gamma)$ we write $(H, H', P) \leq (\tilde{H}, \tilde{H}', \tilde{P})$ if

$\mathbf{c}(H, H', P) < \mathbf{c}(\tilde{H}, \tilde{H}', \tilde{P})$

or if $\mathbf{c}(H, H', P) = \mathbf{c}(\tilde{H}, \tilde{H}', \tilde{P})$ and $|H'/H| > |\tilde{H}'/\tilde{H}|$

or if $\mathbf{c}(H, H', P) = \mathbf{c}(\tilde{H}, \tilde{H}', \tilde{P})$, $|H'/H| = |\tilde{H}'/\tilde{H}|$ and $|H'| \leq |\tilde{H}'|$ (or equivalently $|H| \leq |\tilde{H}|$).

Assume now that $(H, H', P) \leq (\tilde{H}, \tilde{H}', \tilde{P}) \leq (H, H', P)$. We have $\mathbf{c}(H, H', P) = \mathbf{c}(\tilde{H}, \tilde{H}', \tilde{P})$, $|H'/H| = |\tilde{H}'/\tilde{H}|$ and $|H'| = |\tilde{H}'|$ (or equivalently $|H| = |\tilde{H}|$). We can check (in each case) that $H' = \tilde{H}'$, $H = \tilde{H}$, $P = \tilde{P}$. It follows that

(b) $\leq$ is a total order on $B'(S_n)$.

For $\Gamma$ as in (a) we define a map

(c) $B'(\Gamma) \to \mathbf{C}[M(\Gamma)]$

by $(H, H', P) \mapsto \mathbf{s}_{H,H';\Gamma} P$. The image of this map is a subset $B(\Gamma)$ of $\mathbf{C}[M(\Gamma)]$.

**3.2.** For $n \in \{1,2,3,4,5\}$ we define a tableau $\Sigma_n$ with rows (indexed by the various $(H,H',P) \in B'(S_n)$ and ordered according to the total order 3.1(b) of $B'(S_n)$) of the form

$$(?).....\widehat{(x,\sigma)} = \mathbf{s}_{H,H';S_n}P.....\mathbf{s}_{H,H';S_n}P^! = \widehat{(x',\sigma')}$$

where

$(?) = \mathbf{c}(H,H',P)$,

$\widehat{(x,\sigma)}$ (with $(x,\sigma) \in M(S_n)$), is another notation for $\mathbf{s}_{H,H';S_n}P \in B(S_n)$, said to be the left entry of the row,

$\widehat{(x',\sigma')}$ (with $(x',\sigma') \in M(S_n)$, said to be the right entry of the row, is the left entry in the row indexed by $(H,H',P^!) \in B'(S_n)$ with $P^!$ as in 2.4.

(In particular, the tableau associates to each $(H,H',P) \in B'(S_n)$ an element $(x,\sigma) \in M(S_n)$. This association is justified by Theorem 3.4.)

We write $\delta(\widehat{(x,\sigma)}) = N$ if $\widehat{(x,\sigma)}$ appears as the left entry in the $N$-th row of $\Sigma_n$.

We write $\delta^!(\widehat{(x',\sigma')}) = N$ if $\widehat{(x',\sigma')}$ appears as the right entry in the $N$-th row of $\Sigma_n$.

The tableau $\Sigma_1$ is:

$(1).....\widehat{(1,1)} = \mathbf{s}_{1,1;S_1}P_1.....\mathbf{s}_{1,1;S_1}P_1 = \widehat{(1,1)}$

The tableau $\Sigma_2$ is:

$(1).....\widehat{(1,1)} = \mathbf{s}_{1,S_2;S_2}P_{(-1)^0}.....\mathbf{s}_{1,S_2;S_2}P_{-1} = \widehat{(g_2,\epsilon)}$
$(1).....\widehat{(1,\epsilon)} = \mathbf{s}_{1,1;S_2}P_1.....\mathbf{s}_{1,1;S_2}P_1 = \widehat{(1,\epsilon)}$
$(1).....\widehat{(g_2,1)} = \mathbf{s}_{S_2,S_2;S_2}P_1.....\mathbf{s}_{S_2,S_2;S_2}P_1 = \widehat{(g_2,1)}$
$(-1).....\widehat{(g_2,\epsilon)} = \mathbf{s}_{1,S_2;S_2}P_{-1}.....\mathbf{s}_{1,S_2;S_2}P_{(-1)^0} = \widehat{(1,1)}$

The tableau $\Sigma_3$ is:

$(1).....\widehat{(1,1)} = \mathbf{s}_{1,S_3;S_3}P_{\theta^0}.....\mathbf{s}_{1,S_3;S_3}P_{\theta^2} == \widehat{(g_3,\theta^2)}$
$(1).....\widehat{(1,r)} = \mathbf{s}_{1,S_2;S_3}P_{(-1)^0}.....\mathbf{s}_{1,S_2;S_3}P_{-1} = \widehat{(g_2,\epsilon)}$
$(1).....\widehat{(1,\epsilon)} = \mathbf{s}_{1,1;S_3}P_1.....\mathbf{s}_{1,1;S_3}P_1 = \widehat{(1,\epsilon)}$
$(1).....\widehat{(g_2,1)} = \mathbf{s}_{S_2,S_2;S_3}P_1.....\mathbf{s}_{S_2,S_2;S_3}P_1 = \widehat{(g_2,1)}$
$(1).....\widehat{(g_3,1)} = \mathbf{s}_{S_3,S_3;S_3}P_1.....\mathbf{s}_{S_3,S_3;S_3}P_1 = \widehat{(g_3,1)}$
$(-1).....\widehat{(g_2,\epsilon)} = \mathbf{s}_{1,S_2;S_3}P_{-1}.....\mathbf{s}_{1,S_2;S_3}P_{(-1)^0} = \widehat{(1,r)}$
$(\theta).....\widehat{(g_3,\theta)} = \mathbf{s}_{1,S_3;S_3}P_{\theta}.....\mathbf{s}_{1,S_3;S_3}P_{(\theta)^0} = \widehat{(1,1)}$
$(\theta^2).....\widehat{(g_3,\theta^2)} = \mathbf{s}_{1,S_3;S_3}P_{\theta^2}.....\mathbf{s}_{1,S_3;S_3}P_{\theta} = \widehat{(g_3,\theta)}$

The tableau $\Sigma_4$ is:

$(1).....\widehat{(1,1)}=\mathbf{s}_{1,S_4;S_4}P_{i^0}.....\mathbf{s}_{1,S_4;S_4}P_{i^3}=\widehat{(g_4,i^3)}$
$(1).....\widehat{(1,\lambda^1)}=\mathbf{s}_{1,S_3;S_4}P_{\theta^0}.....\mathbf{s}_{1,S_3;S_4}P_{\theta^2}=\widehat{(g_3,\theta^2)}$
$(1).....\widehat{(1,\sigma)}=\mathbf{s}_{1,S_2\times S_2;S_4}P_{1,1}.....\mathbf{s}_{1,S_2\times S_2;S_4}P_{-1,-1}=\widehat{(g'_2,\epsilon)}$
$(1).....\widehat{(1,\lambda^2)}=\mathbf{s}_{1,S_2;S_4}P_{(-1)^0}.....\mathbf{s}_{1,S_2;S_4}P_{-1}=\widehat{(g_2,\epsilon)}$
$(1).....\widehat{(g_2,1)}=\mathbf{s}_{S_2,S_2\times S_2;S_4}P_{(-1)^0}.....\mathbf{s}_{S_2,S_2\times S_2;S_4}P_{-1}=\widehat{(g'_2,r)}$
$(1).....\widehat{(g'_2,1)}=\mathbf{s}_{S_2\times S_2,\Delta;S_4}P_{(-1)^0}.....\mathbf{s}_{S_2\times S_2,\Delta;S_4}P_{-1}=\widehat{(g_4,-1)}$
$(1).....\widehat{(1,\lambda^3)}=\mathbf{s}_{1,1;S_4}P_1.....\mathbf{s}_{1,1;S_4}P_1=\widehat{(1,\lambda^3)}$
$(1).....\widehat{(g_2,\epsilon'')}=\mathbf{s}_{S_2,S_2;S_4}P_1.....\mathbf{s}_{S_2,S_2;S_4}P_1=\widehat{(g_2,\epsilon'')}$
$(1).....\widehat{(g'_2,\epsilon'')}=\mathbf{s}_{S_2\times S_2,S_2\times S_2;S_4}P_1.....\mathbf{s}_{S_2\times S_2,S_2\times S_2;S_4}P_1=\widehat{(g'_2,\epsilon'')}$
$(1).....\widehat{(g_3,1)}=\mathbf{s}_{S_3,S_3;S_4}P_1.....\mathbf{s}_{S_3,S_3;S_4}P_1=\widehat{(g_3,1)}$
$(1).....\widehat{(g'_2,\epsilon')}=\mathbf{s}_{\Delta,\Delta;S_4}P_1.....\mathbf{s}_{\Delta,\Delta;S_4}P_1=\widehat{(g'_2,\epsilon')}$
$(1).....\widehat{(g_4,1)}=\mathbf{s}_{S_4,S_4;S_4}P_1.....\mathbf{s}_{S_4,S_4;S_4}P_1=\widehat{(g_4,1)}$
$(-1).....\widehat{(g_2,\epsilon')}=\mathbf{s}_{1,S_2\times S_2;S_4}P_*.....\mathbf{s}_{1,S_2\times S_2;S_4}P_*=\widehat{(g_2,\epsilon')}$
$(-1).....\widehat{(g_2,\epsilon)}=\mathbf{s}_{1,S_2;S_4}P_{-1}.....\mathbf{s}_{1,S_2;S_4}P_{(-1)^0}=\widehat{(1,\lambda^2)}$
$(-1).....\widehat{(g'_2,r)}=\mathbf{s}_{S_2,S_2\times S_2;S_4}P_{-1}.....\mathbf{s}_{S_2,S_2\times S_2;S_4}P_{(-1)^0}=\widehat{(g_2,1)}$
$(-1).....\widehat{(g_4,-1)}=\mathbf{s}_{S_2\times S_2,\Delta;S_4}P_{-1}.....\mathbf{s}_{S_2\times S_2,\Delta;S_4}P_{(-1)^0}=\widehat{(g'_2,1)}$
$(-1').....\widehat{(g'_2,\epsilon)}=\mathbf{s}_{1,S_2\times S_2;S_4}P_{-1,-1}.....\mathbf{s}_{1,S_2\times S_2;S_4}P_{1,1}=\widehat{(1,\sigma)}$
$(\theta).....\widehat{(g_3,\theta)}=\mathbf{s}_{1,S_3;S_4}P_{\theta}.....\mathbf{s}_{1,S_3;S_4}P_{\theta^0}=\widehat{(1,\lambda^1)}$
$(\theta^2).....\widehat{(g_3,\theta^2)}=\mathbf{s}_{1,S_3;S_4}P_{\theta^2}.....\mathbf{s}_{1,S_3;S_4}P_{\theta}=\widehat{(g_3,\theta)}$
$(i).....\widehat{(g_4,i)}=\mathbf{s}_{1,S_4;S_4}P_i.....\mathbf{s}_{1,S_4;S_4}P_{i^0}=\widehat{(1,1)}$
$(i^3).....\widehat{(g_4,i^3)}=\mathbf{s}_{1,S_4;S_4}P_{i^3}.....\mathbf{s}_{1,S_4;S_4}P_i=\widehat{(g_4,i)}$

The tableau $\Sigma_5$ is:

$(1).....\widehat{(1,1)}=\mathbf{s}_{1,S_5;S_5}P_{\zeta^0}.....\mathbf{s}_{1,S_5;S_5}P_{\zeta^4}=\widehat{(g_5,\zeta^4)}$
$(1).....\widehat{(1,\lambda^1)}=\mathbf{s}_{1,S_4;S_5}P_{i^0}.....\mathbf{s}_{1,S_4;S_5}P_{i^3}=\widehat{(g_4,i^3)}$
$(1).....\widehat{(1,\nu)}=\mathbf{s}_{1,S_3\times S_2;S_5}P_{(-\theta)^0}.....\mathbf{s}_{1,S_3\times S_2;S_5}P_{(-\theta)^5}=\widehat{(g_6,-\theta^2)}$
$(1).....\widehat{(1,\lambda^2)}=\mathbf{s}_{1,S_3;S_5}P_{\theta^0}.....\mathbf{s}_{1,S_3;S_5}P_{\theta^2}=\widehat{(g_3,\epsilon\theta^2)}$
$(1).....\widehat{(g_2,1)}=\mathbf{s}_{S_2,S_3\times S_2;S_5}P_{\theta^0}.....\mathbf{s}_{S_2,S_3\times S_2;S_5}P_{\theta^2}=\widehat{(g_6,\theta^2)}$
$(1).....\widehat{(1,\nu')}=\mathbf{s}_{1,S_2\times S_2;S_5}P_{1,1}.....\mathbf{s}_{1,S_2\times S_2;S_5}P_{-1,-1}=\widehat{(g'_2,\epsilon)}$
$(1).....\widehat{(1,\lambda^3)}=\mathbf{s}_{1,S_2;S_5}P_{(-1)^0}.....\mathbf{s}_{1,S_2;S_5}P_{-1}=\widehat{(g_2,-\epsilon)}$
$(1).....\widehat{(g_2,r)}=\mathbf{s}_{S_2,S_2\times S_2;S_5}P_{(-1)^0}.....\mathbf{s}_{S_2,S_2\times S_2;S_5}P_{-1}=\widehat{(g'_2,r)}$
$(1).....\widehat{(g'_2,1)}=\mathbf{s}_{S_2\times S_2,\Delta;S_5}P_{(-1)^0}.....\mathbf{s}_{S_2\times S_2,\Delta;S_5}P_{-1}=\widehat{(g_4,-1)}$
$(1).....\widehat{(g_3,1)}=\mathbf{s}_{S_3,S_3\times S_2;S_5}P_{(-1)^0}.....\mathbf{s}_{S_3,S_3\times S_2;S_5}P_{-1}=\widehat{(g_6,-1)}$
$(1).....\widehat{(1,\lambda^4)}=\mathbf{s}_{1,1;S_5}P_1.....\mathbf{s}_{1,1;S_5}P_1=\widehat{(1,\lambda^4)}$
$(1).....\widehat{(g_2,\epsilon)}=\mathbf{s}_{S_2,S_2;S_5}P_1.....\mathbf{s}_{S_2,S_2;S_5}P_1=\widehat{(g_2,\epsilon)}$

$(1).....\widehat{(g_2',\epsilon'')} = \mathbf{s}_{S_2\times S_2,S_2\times S_2;S_5}P_1.....\mathbf{s}_{S_2\times S_2,S_2\times S_2;S_5}P_1 = \widehat{(g_2',\epsilon'')}$
$(1).....\widehat{(g_3,\epsilon)} = \mathbf{s}_{S_3,S_3;S_5}P_1.....\mathbf{s}_{S_3,S_3;S_5}P_1 = \widehat{(g_3,\epsilon)}$
$(1).....\widehat{(g_2',\epsilon')} = \mathbf{s}_{\Delta,\Delta;S_5}P_1.....\mathbf{s}_{\Delta,\Delta;S_5}P_1 = \widehat{(g_2',\epsilon')}$
$(1).....\widehat{(g_6,1)} = \mathbf{s}_{S_3\times S_2,S_3\times S_2;S_5}P_1.....\mathbf{s}_{S_3\times S_2,S_3\times S_2;S_5}P_1 = \widehat{(g_6,1)}$
$(1).....\widehat{(g_4,1)} = \mathbf{s}_{S_4,S_4;S_5}P_1.....\mathbf{s}_{S_4,S_4;S_5}P_1 = \widehat{(g_4,1)}$
$(1).....\widehat{(g_5,1)} = \mathbf{s}_{S_5,S_5;S_5}P_1.....\mathbf{s}_{S_5,S_5;S_5}P_1 = \widehat{(g_5,1)}$
$(-1).....\widehat{(g_2,-1)} = \mathbf{s}_{1,S_3\times S_2;S_5}P_{(-\theta)^3}.....\mathbf{s}_{1,S_3\times S_2;S_5}P_{(-\theta)^2} = \widehat{(g_3,\theta^2)}$
$(-1).....\widehat{(g_2,-r)} = \mathbf{s}_{1,S_2\times S_2;S_5}P_*.....\mathbf{s}_{1,S_2\times S_2;S_5}P_* = \widehat{(g_2,-r)}$
$(-1).....\widehat{(g_2,-\epsilon)} = \mathbf{s}_{1,S_2;S_5}P_{-1}.....\mathbf{s}_{1,S_2;S_5}P_{(-1)^0} = \widehat{(1,\lambda^3)}$
$(-1).....\widehat{(g_2',r)} = \mathbf{s}_{S_2,S_2\times S_2;S_5}P_{-1}.....\mathbf{s}_{S_2,S_2\times S_2;S_5}P_{(-1)^0} = \widehat{(g_2,r)}$
$(-1).....\widehat{(g_4,-1)} = \mathbf{s}_{S_2\times S_2,\Delta;S_5}P_{-1}.....\mathbf{s}_{S_2\times S_2,\Delta;S_5}P_{(-1)^0} = \widehat{(g_2',1)}$
$(-1).....\widehat{(g_6,-1)} = \mathbf{s}_{S_3,S_3\times S_2;S_5}P_{-1}.....\mathbf{s}_{S_3,S_3\times S_2;S_5}P_{(-1)^0} = \widehat{(g_3,1)}$
$(-1').....\widehat{(g_2',\epsilon)} = \mathbf{s}_{1,S_2\times S_2;S_5}P_{-1,-1}.....\mathbf{s}_{1,S_2\times S_2;S_5}P_{1,1} = \widehat{(1,\nu')}$
$(\theta).....\widehat{(g_3,\theta)} = \mathbf{s}_{1,S_3\times S_2;S_5}P_{(-\theta)^4}.....\mathbf{s}_{1,S_3\times S_2;S_5}P_{(-\theta)^3} = \widehat{(g_2,-1)}$
$(\theta).....\widehat{(g_3,\epsilon\theta)} = \mathbf{s}_{1,S_3;S_5}P_{\theta}.....\mathbf{s}_{1,S_3;S_5}P_{\theta^0} = \widehat{(1,\lambda^2)}$
$(\theta).....\widehat{(g_6,\theta)} = \mathbf{s}_{S_2,S_3\times S_2;S_5}P_{\theta}.....\mathbf{s}_{S_2,S_3\times S_2;S_5}P_{\theta^0} = \widehat{(g_2,1)}$
$(\theta^2).....\widehat{(g_3,\theta^2)} = \mathbf{s}_{1,S_3\times S_2;S_5}P_{(-\theta)^2}.....\mathbf{s}_{1,S_3\times S_2;S_5}P_{-\theta} = \widehat{(g_6,-\theta)}$
$(\theta^2).....\widehat{(g_3,\epsilon\theta^2)} = \mathbf{s}_{1,S_3;S_5}P_{\theta^2}.....\mathbf{s}_{1,S_3;S_5}P_{\theta} = \widehat{(g_3,\epsilon\theta)}$
$(\theta^2).....\widehat{(g_6,\theta^2)} = \mathbf{s}_{S_2,S_3\times S_2;S_5}P_{\theta^2}.....\mathbf{s}_{S_2,S_3\times S_2;S_5}P_{\theta} = \widehat{(g_6,\theta)}$
$(i).....\widehat{(g_4,i)} = \mathbf{s}_{1,S_4;S_5}P_i.....\mathbf{s}_{1,S_4;S_5}P_{i^0} = \widehat{(1,\lambda^1)}$
$(i^3).....\widehat{(g_4,i^3)} = \mathbf{s}_{1,S_4;S_5}P_{i^3}.....\mathbf{s}_{1,S_4;S_5}P_i = \widehat{(g_4,i)}$
$(\zeta).....\widehat{(g_5,\zeta)} = \mathbf{s}_{1,S_5;S_5}P_{\zeta}.....\mathbf{s}_{1,S_5;S_5}P_{\zeta^0} = \widehat{(1,1)}$
$(\zeta^2).....\widehat{(g_5,\zeta^2)} = \mathbf{s}_{1,S_5;S_5}P_{\zeta^2}.....\mathbf{s}_{1,S_5;S_5}P_{\zeta^3} = \widehat{(g_5,\zeta^3)}$
$(\zeta^3).....\widehat{(g_5,\zeta^3)} = \mathbf{s}_{1,S_5;S_5}P_{\zeta^3}.....\mathbf{s}_{1,S_5;S_5}P_{\zeta^2} = \widehat{(g_5,\zeta^2)}$
$(\zeta^4).....\widehat{(g_5,\zeta^4)} = \mathbf{s}_{1,S_5;S_5}P_{\zeta^4}.....\mathbf{s}_{1,S_5;S_5}P_{\zeta} = \widehat{(g_5,\zeta)}$
$(-\theta).....\widehat{(g_6,-\theta)} = \mathbf{s}_{1,S_3\times S_2;S_5}P_{-\theta}.....\mathbf{s}_{1,S_3\times S_2;S_5}P_{(-\theta)^0} = \widehat{(1,\nu)}$
$(-\theta^2).....\widehat{(g_6,-\theta^2)} = \mathbf{s}_{1,S_3\times S_2;S_5}P_{(-\theta)^5}.....\mathbf{s}_{1,S_3\times S_2;S_5}P_{(-\theta)^4} = \widehat{(g_3,\theta)}$

**3.3.** Let $n\in\{1,2,3,4,5\}$. For $(x,s)\in M(S_n)$ we write

$$\widehat{(x,\sigma)} = \sum_{(y,\tau)\in M(S_n)} \gamma_{(x,\sigma),(y,\tau)}(y,\tau)$$

where $\gamma_{(x,\sigma),(y,\tau)}\in\mathbf{C}$.

**Theorem 3.4.** *If $\gamma_{(x,\sigma),(y,\tau)} \ne 0$ then $\delta(\widehat{(y,\tau)}) \le \delta(\widehat{(x,\sigma)})$. We have $\gamma_{(x,\sigma),(y,\tau)} \in$ $\mathbf{N}$ except for a single pair $(x,\sigma),(y,\tau)$ (for $n=5$) when $\gamma_{(x,\sigma),(y,\tau)} = -(\zeta^2+\zeta^3) \in$ $\mathbf{R}_{\ge 0}$. We have also $\gamma_{(x,\sigma),(x,\sigma)} = 1$.*

In the following formulas, $*$ represents an $\mathbf{N}$-linear combination of elements $(x,\sigma) \in M(S_n)$ corresponding to rows in $\Sigma_n$ that start with (1). (Note that $*$ is easily computable in each case.)

In $\Sigma_3$ we have
$\widehat{(g_3,\theta)} = (g_3,\theta) + *$
$\widehat{(g_3,\theta^2)} = (g_3,\theta) + (g_3,\theta^2) + *$.
In $\Sigma_4$ we have
$\widehat{(g_3,\theta)} = (g_3,\theta) + *$
$\widehat{(g_3,\theta^2)} = (g_3,\theta) + (g_3,\theta^2) + *$
$\widehat{(g_4,i)} = (g_4,i) + *$
$\widehat{(g_4,i^3)} = (g_4,i) + (g_4,i^3) + *$.
In $\Sigma_5$ we have
$\widehat{(g_3,\theta)} = (g_3,\theta) + *$
$\widehat{(g_3,\epsilon\theta)} = (g_3,\epsilon\theta) + (g_3,\theta) + *$
$\widehat{(g_6,\theta)} = (g_6,\theta) + (g_3,\theta) + *$
$\widehat{(g_3,\theta^2)} = (g_3,\theta) + (g_3,\theta^2) + *$
$\widehat{(g_3,\epsilon\theta^2)} = (g_3,\epsilon\theta) + (g_3,\theta) + (g_3,\epsilon\theta^2) + (g_3,\theta^2) + *$
$\widehat{(g_6,\theta^2)} = (g_6,\theta) + (g_6,\theta^2) + (g_3,\theta) + (g_3,\theta^2) + *$
$\widehat{(g_4,i)} = (g_4,i) + *$
$\widehat{(g_4,i^3)} = (g_4,i) + (g_4,i^3) + *$
$\widehat{(g_5,\zeta} = (g_5,\zeta) + *$
$\widehat{(g_5,\zeta^2} = -(\zeta^2+\zeta^3)(g_5,\zeta) + (g_5,\zeta^2) + *$
$\widehat{(g_5,\zeta^3} = (g_5,\zeta^2) + (g_5,\zeta^3) + *$
$\widehat{(g_5,\zeta^4} = (g_5,\zeta) + (g_5,\zeta^2) + (g_5,\zeta^3) + (g_5,\zeta^4) + *$
$\widehat{(g_6,-\theta)} = (g_6,-\theta) + *$
$\widehat{(g_6,-\theta^2)} = (g_6,-\theta) + (g_6,-\theta^2) + *$.

All other elements $\widehat{(x,\sigma)}$ in $M(S_n)$ are those listed in [L20]. From this we see that the theorem holds.

**Corollary 3.5.** *$B_{(}S_n)$ is a $\mathbf{C}$-basis of $\mathbf{C}[M(S_n)]$ related to the basis formed by $M(S_n)$ by a triangular matrix (for the total order $\le$ on $B'(S_n)$) with $1$ on diagonal and with entries in $\mathbf{N}$ (except for one entry $-(\zeta^2+\zeta^3)$ in the case $n=5$). In particular the map $M(S_n) \to B(S_n)$, $(x,\sigma) \mapsto \widehat{(x,\sigma)}$ is a bijection.*

**Corollary 3.6.** *There is a unique bijection $\epsilon : M(S_n) \to B(S_n)$ such that for any $(x,\sigma) \in M(S_n)$, $(x,\sigma)$ appears with nonzero coefficient in $\epsilon(x,\sigma)$. Moreover that coefficient is* 1.

**Corollary 3.7.** *The map $B'(S_n) \to B(S_n)$ defined by 3.1(c) is a bijection.*

**3.8.** We see that $B(S_n)$ has the properties 0.2(II),(III). Using 0.1(b) and 2.3 we see that it also has property 0.2(I).

**3.9.** For $b = \mathbf{s}_{H,H';S_n}P \in B(S_n)$ we set $b^! = \mathbf{s}_{H,H'}P^! \in B(S_n)$. Now $b \mapsto b^!$ is a bijection $B(S_n) \to B(S_n)$.

For $b \in B(S_n)$ we can write uniquely

$$A_{S_n}(b^!) = \sum_{b' \in B(S_n)} f_{b,b'} b'^!.$$

where $f_{b,b'} \in \mathbf{C}$. For $b = \mathbf{s}_{H,H';S_n}P$, $(H,H',P) \in B'(S_n)$ we set $\mathbf{c}(b) = \mathbf{c}_{H,H',P}$.

**Theorem 3.10.** *(a) If $f_{b,b'} \neq 0$ then $\mathbf{c}(b'^!) \le \mathbf{c}(b^!)$. If in addition we have $\mathbf{c}(b'^!) = \mathbf{c}(b^!)$ then $b' = b$.*

*(b) For $b \in B(S_n)$ we set $f_{b,b} = sg(b)$. We have $sg(b) \in \{\pm 1\}$.*

Note that (b) is a consequence of (a) since it is known [L79] that $A_{S_n}^2 = 1$.

We now prove (a).

Assume first that $n = 5$ and $b = \widehat{(1,1)} = (1,1)$. From 2.3(a) we see that $A_{S_5}(b)$ is a linear combination of terms of the form

(c) $(g_5,\zeta) + (g_5,\zeta^2) + (g_5,\zeta^3) + (g_5,\zeta^4)$,
$(g_6,-\theta) + (g_6,-\theta^2)$, $(g_6,\theta) + (g_6,\theta^2)$,
$(g_3,\theta) + (g_3,\theta^2)$, $(g_3,\epsilon\theta) + (g_3,\epsilon\theta^2)$,
$(g_4,i) + (g_4,i^3)$, $(g_5,1),(g_6,1),(g_6,-1),(g_3,1),(g_3,\epsilon)$,
$(g_4,1),(g_4,-1)$, $(g_2,\rho)$, $(g_2',\rho')$, $(1,\rho'')$
for various $\rho,\rho',\rho''$.

Each of these terms is a linear combination of elements $\widehat{(x,\sigma)}$ which have $\delta^! \le 34$. (We use 3.4.) Thus (a) holds for $b = (1,1)$ (we have $\delta^!(b) = 34$).

Applying $\mathbf{s}_{1,C_5;S_5}$ to 1.4(e) and using 2.2 and 0.1(a), we see that

(d) $A_{S_5}P_{\zeta^4} = P_{\zeta^4}$;

thus (a) holds for $b = P_{\zeta^4}$ (we have $\delta^!(b) = 0$).

Applying $\mathbf{s}_{1,C_5;S_5}$ to 1.4(g) and using 2.2 and 0.1(a), we see that

(e) $A_{S_5}(P_{\zeta^3} + \mathcal{H}_5 - \mathcal{H}_5') - (P_{\zeta^3} + \mathcal{H}_5 - \mathcal{H}_5')$ is a linear combination of $P_{\zeta^4}$ (with $\delta^! = 0$) and of terms $\mathbf{s}_{1,C_5;S_5}\mathcal{T}, \mathbf{s}_{1,C_5;S_5}\mathcal{T}', \mathbf{s}_{1,C_5;S_5}(0,0)$.

Each of these terms is a linear combination of $\widehat{(x,\sigma)}$ with $\delta^! \le 34$. (We use 3.4.) Also $A_{S_5}(\mathcal{H}_5 - \mathcal{H}_5')$ is a sum of terms $A_{S_5}(1,?)$; these are linear combination of terms in (a); each of these terms is a linear combination of elements $\widehat{(x,\sigma)}$ which have $\delta^! \le 34$. (We use 3.4.) Moreover $\mathcal{H}_5 - \mathcal{H}_5'$ is a linear combination of elements $\widehat{(x,\sigma)}$ which have $\delta^! \le 34$. Hence (a) holds for $b = P_{\zeta^3}$. (We have $\delta^!(b) = 35$.)

Applying $\mathbf{s}_{1,C_5;S_5}$ to 1.4(i) and using 2.2 and 0.1(a), we see that

(f) $A_{S_5}P_{\zeta^2}-P_{\zeta^2}$ is a linear combination of $P_{\zeta^4}$ (with $\delta^!=0$) and of terms $\mathbf{s}_{1,C_5;S_5}\mathcal{T},\mathbf{s}_{1,C_5;S_5}\mathcal{T}',\mathbf{s}_{1,C_5;S_5}(0,0)$.

Each of these terms is a linear combination of $\widehat{(x,\sigma)}$ with $\delta^!\le 34$. (We use 3.4.) Hence (a) holds for $b=P_{\zeta^2}$. (We have $\delta^!(b)=36$.)

Applying $\mathbf{s}_{1,C_5;S_5}$ to 1.4(j) and using 2.2 and 0.1(a), we see that

(g) $A_{S_5}P_\zeta+P_\zeta$ is a linear combination of $P_{\zeta^4}$ (with $\delta^!=0$), $P_{\zeta^3}$ (with $\delta^!=35$), $P_{\zeta^2}$ (with $\delta^!=36$) and of terms $\mathbf{s}_{1,C_5;S_5}\mathcal{T},\mathbf{s}_{1,C_5;S_5}\mathcal{T}',\mathbf{s}_{1,C_5;S_5}(0,0)$.

Each of these terms is a linear combination of $\widehat{(x,\sigma)}$ with $\delta\le 34$. (We use 3.4.) Hence (a) holds for $b=P_\zeta$. (We have $\delta^(b)=37$.)

Assume next that $n=4$ and $b=\widehat{(1,1)}=(1,1)$. From 2.3(a) we see that $A_{S_4}(b)$ is a linear combination of terms of the form

(h) $(g_3,\theta)+(g_3,\theta^2)$, $(g_4,i)+(g_4,-i)$, $(g_3,1)$, $(g_4,1),(g_4,-1)$, $(g_2,\rho)$, $(g_2',\rho')$, $(1,\rho'')$ for various $\rho,\rho',\rho''$.

Each of these terms is a linear combination of elements $\widehat{(x,\sigma)}$ which have $\delta^!\le 20$. (We use 3.4.) Thus (a) holds for $b=(1,1)$ (we have $\delta^!(b)=20$).

Applying $\mathbf{s}_{1,C_4;S_4}$ to 1.3(b) and using 2.2 and 0.1(a), we see that

(i) $A_{S_4}P_{i^3}=P_{i^3}$;

thus (a) holds for $b=P_{i^3}$ (we have $\delta^!(b)=0$).

Applying $\mathbf{s}_{1,C_4;S_4}$ to 1.3(e) and using 2.2 and 0.1(a), we see that

(j) $A_{S_4}(P_i+\mathcal{H}_4-\mathcal{H}_4')-(P_i+\mathcal{H}_4-\mathcal{H}_4')$ is a linear combination of $P_{i^3}$ and of terms

$$\mathbf{s}_{1,C_4;S_4}(0,0),\mathbf{s}_{1,C_4;S_4}(0,2),\mathbf{s}_{1,C_4;S_4}(2,0),\mathbf{s}_{1,C_4;S_4}(2,2).$$

Each of these terms is a linear combination of elements $\widehat{(x,\sigma)}$ which have $\delta^!\le 20$. (We use 3.4.) Also $A_{S_4}(\mathcal{H}_4-\mathcal{H}_4')$ is a sum of terms $A_{S_4}(1,?)$; these are linear combination of terms in (h); each of these terms is a linear combination of elements $\widehat{(x,\sigma)}$ which have $\delta^!\le 20$. (We use 3.4.) Moreover $\mathcal{H}_4-\mathcal{H}_4'$ is a linear combination of elements $\widehat{(x,\sigma)}$ which have $\delta^!\le 20$. Hence (a) holds for $b=P_i$ (we have $\delta^!(b)=21$).

The general case can be reduced to the cases already treated using 0.1(a). For example, let $b\in B_n$ be such that $\mathbf{c}(b)=(1)$. We have $b^!=\mathbf{s}_{H,H';S_n}P$ where $A_{H'/H}P=P$. Using 0.1(a) we have

$$A_{S_n}(b^!)=A_{S_n}\mathbf{s}_{H,H';S_n}P=\mathbf{s}_{H,H';S_n}A_{H'/H}P=\mathbf{s}_{H,H';S_n}P=b^!$$

hence (a) is verified in this case. This completes the proof of (a).

We see that 0.2(IV) holds.

**3.11.** We say that $b\in B(S_n)$ is Lagrangian if $A_{S_n}(b)=b$.

(a) *The Lagrangian elements of $B(S_n)$ are precisely the right entries of the rows of $\Sigma_n$ which start with* (1).

## 4. The function sg : $B(S_n) \to \{\pm 1\}$

**4.1.** Let $\Gamma$ be $S_4$ or $S_5$. Let

$$\tau_\Gamma = \sum_{b \in B(\Gamma)} \mathrm{sg}(b) \in \mathbf{Z}.$$

We have

$$\tau_\Gamma = \mathrm{tr}(A_\Gamma : \mathbf{C}[M(\Gamma)] \to \mathbf{C}[M(\Gamma)]) = \sum_{(x,\sigma) \in M(\Gamma)} \{(x,\sigma),(x,\sigma)\}.$$

When $\Gamma = S_4$, the last sum is

$$\tau_\Gamma = 1/24 + 3/8 + 1/6 + 3/8 + 1/24 + 1/2 + 3/8 + 1/2 + 2/3$$
$$+ 1/2 + 3/8 + 3/8 + 1/2 + 1/2 + 1/2 + 1/2 + 3/8 + 2/3 + 2/3 + 1/2 + 1/2 = 9.$$

When $\Gamma = S_5$, the analogous sum is

$$\tau_\Gamma = 1/120 + 2/15 + 5/24 + 3/10 + 5/24 + 2/15 + 1/120$$
$$+ 1/3 + 1/3 + 1/3 + 3/8 + 3/8 + 1/3 + 1/3$$
$$+ 1/3 + 1/2 + 4/5 + 3/8 + 1/3 + 1/3 + 1/2 + 1/2 + 1/3 + 1/3 + 3/8$$
$$+ 1/3 + 1/3 + 1/3 + 1/3 + 1/3 + 1/3 + 1/2 + 1/2$$
$$+ (2 + \zeta + \zeta^{-1})/5 + (2 + \zeta^2 + \zeta^{-2})/5$$
$$+ (2 + \zeta^3 + \zeta^{-3})/5 + (2 + \zeta^4 + \zeta^{-4})/5 + 1/3 + 1/3 = 13.$$

**4.2.** Let $n \in \{1,2,3,4,5\}$. We have the following result.

**Proposition 4.3.** *Let $b \in B(S_n)$. Write $b = \mathbf{s}_{H,H';S_n} P$ where $(H,H',P) \in B'(S_n)$. We have $sg(b) = sg'(P)$ (with $sg'(P)$ as in 2.1).*

If $H'/H$ is one of $S_1, S_2, S_3, S_3 \times S_2, S_2 \times S_2$ one can check from the definitions that the result holds. Thus we can asssume that $H'/H$ is $S_5$ or $S_4$ so that $H = 1$ and $H'$ is $S_5$ or $S_4$. From the definitions one can check that $\mathrm{sg}(b) = \mathrm{sg}(P)$. It remains to show that $\mathrm{sg}(P) = \mathrm{sg}'(P)$.

If $n = 4$ and $P \in Prim(S_4)$, $P \neq (1,1)$ we have $\mathrm{sg}(P) = \mathrm{sg}'(P)$, by 3.10(g),(h). Then in the sum $\sum_{b \in B(S_4)} \mathrm{sg}(b) = 9$ (see 4.1) all terms are known except $\mathrm{sg}(1,1)$ which is then necessarily $1 = \mathrm{sg}'(1,1)$.

If $n = 5$ and $P \in Prim(S_5)$, $P \neq (1,1)$ we have $\mathrm{sg}(P) = \mathrm{sg}'(P)$, by 3.10(c),(d),(e),(f). Then in the sum $\sum_{b \in B(S_5)} \mathrm{sg}(b) = 13$ (see 4.1) all terms are known except $\mathrm{sg}(1,1)$ which is then necessarily $-1 = \mathrm{sg}'(1,1)$. This proves the proposition.

## 5. Relation with unipotent representations

**5.1.** In this section we assume that $G$ is a simple algebraic group over an algebraic closure of the finite field $F_q$ with $q$ elements and that $G$ is provided with a split $F_q$-rational structure. Let $G(F_q)$ be the finite group of $F_q$-rational points of $G$. Let $W$ be the Weyl group of $W$. Then $\mathrm{Irr}(W)$ is a disjoint union of subsets called families, see [L79]. Let $\mathcal{U}$ be the subset of $\mathrm{Irr}(G(F_q))$ consisting of all unipotent representations of $G(F_q)$. It is known [L84] that $\mathcal{U} = \sqcup_c \mathcal{U}_c$ where $c$ runs over the families of $W$. To each family $c$, one can associate a finite group $\Gamma_c$ so that $\mathcal{U}_c$ is in bijection

(a) $\rho_{x,\sigma} \leftrightarrow (x,\sigma)$

with $M(\Gamma_c)$, see [L84].

**5.2.** In the remainder of this section we assume that $G$ is of type $G_2, F_4$ (resp. $E_8$) and $c$ is the unique family such that $\Gamma_c = S_n$ with $n = 3, 4$ (resp. 5).

For $(x,\sigma),(y,\tau)$ in $M(\Gamma_c)$ let $\gamma_{(x,\sigma),(y,\tau)} \in \mathbf{C}$ be as in 3.3. By 3.4 we have $\gamma_{(x,\sigma),(y,\tau)} \in \mathbf{R}_{\geq 0}$. Let

$$\hat{\rho}_{x,\sigma} = \sum_{(y,\tau)\in M(\Gamma_c)} \gamma_{(x,\sigma),(y,\tau)} \rho_{y,\tau}$$

an element of the Grothendieck group of $G(F_q)$-modules tensored with $\mathbf{R}$.

**Theorem 5.3.** $\dim \hat{\rho}_{x,\sigma}$ *is a polynomial in $q$ with coefficients in $\mathbf{R}_{\geq 0}$. (Here we regard $q$ as a variable.)*

According to [DL76,7.6, L84,4.23], for any $(x,\sigma) \in M(\Gamma_c)$ we have

$$\text{(a)} \qquad \dim \rho_{x,\sigma} = \sum_{E\in c, i\in\mathbf{N}} < m_E, (x,\sigma) > (E : \bar{S}_i)_W q^i$$

where $E \mapsto m_E$ is an imbedding $c \to M_{\Gamma_c}$ and $\oplus_i \bar{S}_i$ is a certain $\mathbf{N}$-graded finite dimensional $W$-module associated to $W$. Using (a) we have

$$\dim \hat{\rho}_{x,\sigma} = \sum_{(y,\tau)\in M(G_c)} \gamma_{(x,\sigma),(y,\tau)} \sum_{E\in c, i\in\mathbf{N}} < m_E, (y,\tau) > (E : \bar{S}_i)_W q^i$$

$$\text{(b)} \qquad = \sum_{E\in c, i\in\mathbf{N}} \gamma'_E (E : \bar{S}_i)_W q^i$$

where

$$\gamma'_E = \sum_{(y,\tau)\in M(\Gamma_c)} < m_E, (y,\tau) > g_{(x,\sigma),(y,\tau)}.$$

From the bipositivity of $\widehat{(x,\sigma)}$ we see that $\gamma'_E \in \mathbf{R}_{\geq 0}$ for any $E \in c$. The theorem follows.

**5.4.** It is known that there is a unique $a(c) \in \mathbf{N}$ such that for any $E \in c$ we have:
$(E : \bar{S}_i) = 0$ for any $i < a(c)$,,
$(E : \bar{S}_{a(c)}) = 0$ if $m_E \neq (1,1)$,
$(E : \bar{S}_{a(c)}) = 1$ if $m_E = (1,1)$.

**Theorem 5.5.** *For any $(x,\sigma) \in M(S_n)$, $\dim \hat{\rho}_{x,\sigma}$ is of the form $e_{x,\sigma} q^{a(c)}$+higher powers of $q$ where*

$$e_{x,\sigma} = \sum_{(y,\tau) \in M(\Gamma_c)} < (1,1),(y,\tau) > \gamma_{(x,\sigma),(y,\tau)} \in \mathbf{R}_{\geq 0}.$$

This follows immediately from 5.3(b) using results in 5.4.

**Theorem 5.6.** *Let $(x,\sigma) \in \Gamma_c$. We have $\widehat{(x,\sigma)} = \mathbf{s}_{H,H';\Gamma_c} P$ where $(H,H') \in Z(S_n), P \in Prim(H'/H)$. If $\mathbf{c}(H,H',P) = (1)$ then $e_{x,\sigma} = |H|/|H'|$ and*

$$\dim \hat{\rho}_{x,\sigma} \in (|H|/|H'|)\mathbf{N}[q].$$

This is proved by direct verification in each case.

**5.8.** Assume that $G$ is of type $G_2$ so that $\Gamma_c = S_3$. For each $(x,\sigma) \in M(S_3)$ we give below the value of $\dim \hat{\rho}_{x,\sigma}$.
$(1,1)...(q/6)(q+1)^2(q^2+q+1)$
$(1,r)...(q/2)(q^2+1)(q^2+q+1)$
$(g_3,\theta)...(q/2)(q+1)(q^3+1)$
$(g_2,1)...q(q^4+q^3+q^2+q+1)$
$(g_3,1)...q(q^4+q^3+q^2+q+1)$
$(1,\epsilon)...q(q^2+1)^2$
$(g_2,\epsilon)...q(q^4+q^2+1)$
$(g_3,\theta^2)...q(q^4+1)$.

**5.9.** We now give some examples of values of $e_{x,\sigma}$ when $\widehat{(x,\sigma)}$ is primitive.
$(n=3)$: $e_{g_3,\theta^2} = 1, e_{g_3,\theta} = 1/2, e_{1,1} = 1/6$.
$(n=4)$: $e_{g_4,i^3} = 1, e_{g_4,i} = 1/2, e_{1,1} = 1/24$.
$(n=5)$: $e_{g_5,\zeta^4} = 1, e_{g_5,\zeta^3} = 1/2, e_{1,1} = 1/120$.

**5.10.** Let $\mathcal{U}_c^{cusp}$ be the set of (unipotent) cuspidal representations in $\mathcal{U}^c$. Let $\eta$ be the composition

$$\mathcal{U}_c^{cusp} \to \mathcal{U}_c \to M(\Gamma_c) = M(S_n) \to B(S_n) \to B'(S_n)$$

(the first map is the inclusion, the second map is the bijection 5.1(a), the third map is the bijection $\epsilon$ in 3.6, the fourth map is the inverse of the bijection 3.7). One can check that:

(a) *If $\rho \in \mathcal{U}_c^{cusp}$ then $\eta(\rho)$ is of the form $(H,H',P)$ where $H=1$. Moreover $\rho \mapsto \mathbf{c}(H,H',P)$ is bijection $\mathcal{U}_c^{cusp} \xrightarrow{\sim} Y_n$.*

We describe this bijection by a table in each case. In each case the rows of the table are of the form
$(\xi).....(x,\sigma).....(1,H',P)$
where $(\xi) \in Y_n, (x,\sigma) \in M(S_n), (1,H',P) \in B'(S_n)$.

Table for $n = 3$

$(1).....(1, \epsilon).....(1, 1, P_1)$
$(-1).....(g_2, \epsilon).....(1, S_2, P_{-1})$
$(\theta).....(g_3, \theta).....(1, S_3, P_\theta)$
$(\theta^2).....(g_3, \theta^2).....(1, S_3, P_{\theta^2})$

Table for $n = 4$

$(1).....(1, \lambda^3).....(1, 1, P_1)$
$(-1).....(g_2, \epsilon).....(1, S_2, P_{-1})$
$(-1').....(g_2', \epsilon).....(1, S_2 \times S_2, P_{-1,-1})$
$(\theta).....(g_3, \theta).....(1, S_3, P_\theta)$
$(\theta^2).....(g_3, \theta^2).....(1, S_3, P_{\theta^2})$
$(i).....(g_4, i).....(1, S_4, P_i)$
$(i^3).....(g_4, i^3).....(1, S_3, P_{i^3})$

Table for $n = 5$

$(1).....(1, \lambda^4).....(1, 1, P_1)$
$(-1).....(g_2, -\epsilon).....(1, S_2, P_{-1})$
$(-1').....(g_2', \epsilon).....(1, S_2 \times S_2, P_{-1,-1})$
$(\theta).....(g_3, \epsilon\theta).....(1, S_3, P_\theta)$
$(\theta^2).....(g_3, \epsilon\theta^2).....(1, S_3, P_{\theta^2})$
$(i).....(g_4, i).....(1, S_4, P_i)$
$(i^3).....(g_4, i^3).....(1, S_4, P_{i^3})$
$(\zeta).....(g_5, \zeta).....(1, S_5, P_\zeta)$
$(\zeta^2).....(g_5, \zeta^2).....(1, S_5, P_{\zeta^2})$
$(\zeta^3).....(g_5, \zeta^3).....(1, S_5, P_{\zeta^3})$
$(\zeta^4).....(g_5, \zeta^4).....(1, S_5, P_{\zeta^4})$
$(-\theta).....(g_6, -\theta).....(1, S_3 \times S_2, P_{-\theta}$
$(-\theta^2).....(g_6, -\theta^2).....(1, S_3 \times S_2, P_{(-\theta)^5}$

Department of Mathematics, M.I.T., Cambridge, MA 02139